\documentclass[a4paper]{amsart}
\usepackage{amsmath,amssymb,amscd}
\usepackage{graphicx}
\usepackage[dvips]{epsfig}
\usepackage{url}
\usepackage{xspace}

\newcommand{\extcurves}{\texttt{extcurves}\xspace}
\newcommand{\CyclePainter}{\texttt{CyclePainter}\xspace}

\newtheorem{theorem}{Theorem}

\allowdisplaybreaks \numberwithin{equation}{section}

\def\Tr{{\rm Tr}}

\newcommand{\dif}{\mathrm{d}}
\newcommand{\mi}{\mathrm{i}}

\title[Klein's curve]%
  {Klein's curve}

\author{H.W. Braden and T.P. Northover}

\address{School of Mathematics, Edinburgh University, Edinburgh.}
\email{hwb@ed.ac.uk}
\address{School of Mathematics, Edinburgh University, Edinburgh.}
\email{T.P.Northover@sms.ed.ac.uk}

\begin{document}

\maketitle

\begin{abstract}
Riemann surfaces with symmetries arise in many studies of
integrable systems. We illustrate new techniques in investigating
such surfaces by means of an example. By giving an homology basis
well adapted to the symmetries of Klein's curve, presented as a
plane curve, we derive a new expression for its period matrix.
This is explicitly related to the hyperbolic model and results of
Rauch and Lewittes.
\end{abstract}

\section{Introduction}

Riemann surfaces with their differing analytic, algebraic,
geometric and topological perspectives have long been been objects
of fascination. In recent years they have appeared in many guises
in mathematical physics (for example string theory and
Seiberg-Witten theory) with integrable systems being a unifying
theme \cite{BK}. In the modern approach to integrable systems a
spectral curve, the Riemann surface in view, plays a fundamental
role: one, for example, constructs solutions to the integrable
system in terms of analytic objects (the Baker-Akhiezer functions)
associated with the curve; the moduli of the curve play the role
of actions for the integrable system. Although a reasonably good
picture exists of the geometry and deformations of these systems
it is far more difficult to construct actual solutions -- largely
due to the transcendental nature of the construction. An actual
solution requires the period matrix of the curve, and there may be
(as in the examples of monopoles, harmonic maps or the AdS-CFT
correspondence) constraints on some of these periods. Fortunately
computational aspects of Riemann surfaces are also developing. The
numerical evaluation of period matrices (suitable for constructing
numerical solutions) is now straightforward but analytic studies
are far more difficult.

When a curve has symmetries the period matrix of the curve is
constrained and this may lead to simplifications, such as allowing
the Baker-Akhiezer functions to be expressed in terms of lower
genus expressions \cite{bbm83}. These simplifications arise (as
will be described below) when an homology basis has been chosen to
reflect the symmetries. For example, when studying magnetic
monopoles with spatial symmetry \cite{hmm95} such a basis
 enables a reduction of the associated Baker-Akhiezer functions
\cite{bren07, bren09, b10}. Finding such an adapted homology basis
is however both a theoretical and practical problem \cite{jg07,
d10}. Here we shall describe a result obtained while developing
new tools to investigate such adapted homology bases. The details
of this software are not the focus here, though we shall point out
where we have employed this in what follows. Rather we believe the
following result both illustrative and of independent interest.

\begin{theorem} The period matrix for Klein's curve (in the homology
basis described below) takes the form \begin{equation} \tau =
   \frac{1}{2}
  \begin{pmatrix}
    e & 1 & 1 \\
    1 & e & 1 \\
    1 & 1 & e
  \end{pmatrix}
  \label{eq:kleinperiod}
\end{equation}
with $e = \frac{-1+\mi\sqrt{7}}{2}$.

\end{theorem}

To put the result in context we note that Levy \cite{l99} has
assembled a very good general introduction to Klein's curve and
various authors have studied the period matrix itself. Since
Riemann the construction of period matrices has been an exercise,
albeit a hard one. The two ingredients, a canonical homology basis
and the integration of a basis of holomorphic differentials along
these, each present different difficulties. Hurwitz
\cite{hur}[p123], Klein and Fricke \cite{fk}[p595], and Baker
\cite{b07} all calculated the periods of a basis of holomorphic
differentials for a set of cycles on Klein's curve. The first
three authors obtained a matrix whose imaginary part was not
positive definite (and so not a period matrix) while Baker
obtained a diagonal matrix (also not a period matrix). Baker
explicitly noted (\cite{b07} footnote p267) that he was not
claiming these cycles formed a canonical homology basis; in fact
they are not. Rauch and Lewittes \cite{rl69}, using a hyperbolic
model, calculated the period matrix to be
\begin{displaymath}
  \tau_{RL} =
  \begin{pmatrix}
    \frac{-1+3\mi\sqrt{7}}{8} &
    \frac{-1-\mi\sqrt{7}}{4} &
    \frac{-3+\mi\sqrt{7}}{8}
    \\
    \frac{-1-\mi\sqrt{7}}{4} &
    \frac{1+\mi\sqrt{7}}{2} &
    \frac{-1-\mi\sqrt{7}}{4}
    \\
    \frac{-3+\mi\sqrt{7}}{8} &
    \frac{-1-\mi\sqrt{7}}{4} &
    \frac{7+3\mi\sqrt{7}}{8}
  \end{pmatrix}.
\end{displaymath}
Given the Rauch-Lewittes result it is a simple matter to establish
the theorem. In terms of the symplectic matrix
\begin{equation} \label{eq:rlsymp}
  M =
  \begin{pmatrix}
    1 & -1 & 0 & 1 & -1 & 0 \\
    0 & -1 & 1 & 0 & -1 & 1 \\
    -1 & -1 & 0 & 0 & 0 & 1 \\
    0 & 1 & 0 & 0 & 0 & -1 \\
    0 & 0 & 0 & 0 & 0 & 1 \\
    -1 & 0 & 0 & -1 & 0 & 0
  \end{pmatrix}=
  \begin{pmatrix}
    A & B \\
    C & D
  \end{pmatrix}
\end{equation}
we find
\[
  \tau_{RL}=(C+D\tau)(A+B\tau)^{-1}.
\]
Thus the period matrices are symplectically related and the
theorem established. Of course the hard work is hidden in the
construction of the symplectic transformation moving from our
homology basis to that of Rauch-Lewittes. In what follows we will
find an homology basis that yields (\ref{eq:kleinperiod}) directly
and then relate this to Rauch and Lewittes'. Other period matrices
for Klein's curve are described later in the paper including the
symmetric result of \cite{rg95} obtained from a hyperbolic model.
We will however begin with a plane curve model as that is how
spectral curves are presented in the integrable system context.

To conclude this introduction we mention the two pieces of
software\footnote{Located at
\url{http://gitorious.org/riemanncycles}.} that have been written
to aid calculations involving intersection numbers of paths on
Riemann surfaces,  \extcurves and \CyclePainter. Such intersection
numbers arise frequently in matters relating to homology and
period matrices and are responsible for the matrices that appear
in this paper. \extcurves accepts piecewise linear descriptions of
homology cycles on a surface (not passing through the branch
points) and calculates their intersection number. In essence this
is just a matter of keeping track of orientations and totalling
the contribution from each intersection, but subtleties arise when
paths don't intersect transversely. It also has peripheral
functions to deal with common applications of the intersection
numbers, such as finding transformations between canonical
homology bases, and to allow interaction with Maple's existing
routines in \texttt{algcurves}. \CyclePainter provides a graphical
interface (resembling a basic paint program) for constructing and
visualising cycles usable in \extcurves, easing the burden of
demanding they are piecewise linear. It provides a 2 dimensional
representation of cycles, with colour corresponding to sheet
information.

\section{Symmetries and the period matrix}
Before turning to the specific setting of Klein's curve let us
recall how symmetries restrict the form of the period matrix. This
will outline our strategy in what follows. Let
$\{\gamma_*\}:=\{\mathfrak{a}_*,\mathfrak{b}_*\}$ be a canonical
homology basis and $\{\omega_j\}$ be a basis of holomorphic
differentials for our Riemann surface $\mathcal{C}$. Set
\begin{equation*}\begin{pmatrix}\oint_{\mathfrak{a}_i}\omega_j\\
\oint_{\mathfrak{b}_i}\omega_j\end{pmatrix}=
\begin{pmatrix}\mathcal{A}\\
\mathcal{B}\end{pmatrix}=\begin{pmatrix}1\\
\tau\end{pmatrix}\mathcal{A}
\end{equation*}
with $\tau=\mathcal{B}\mathcal{A}\sp{-1}$ the
($\mathfrak{a}$-normalized) period matrix. Let $\sigma$ be an
automorphism of $\mathcal{C}$. Then $\sigma$ acts on
$H_1(\mathcal{C},\mathbb{Z})$ and the holomorphic differentials by
\[
  \sigma_\ast\begin{pmatrix}\mathfrak{a}_i\\
  \mathfrak{b}_i\end{pmatrix}=M\begin{pmatrix}\mathfrak{a}_i\\
  \mathfrak{b}_i\end{pmatrix}:=
  \begin{pmatrix}{A}&B\\
    C&D\end{pmatrix}\begin{pmatrix}\mathfrak{a}_i\\
    \mathfrak{b}_i
  \end{pmatrix},
  \qquad \sigma\sp\ast\omega_j=\omega_k L\sp{k}_j,
\]
where $M\in Sp(2g,\mathbb{Z})$. The fundamental
identity
\begin{equation*}\oint_{\sigma_\ast\gamma}\omega=\oint_\gamma\sigma\sp\ast\omega
\end{equation*}
then leads to the constraint
\begin{align}
\begin{pmatrix}{A}&B\\
C&D\end{pmatrix}\begin{pmatrix}\mathcal{A}\\
\mathcal{B}\end{pmatrix}&=\begin{pmatrix}\mathcal{A}\\
\mathcal{B}\end{pmatrix}L \label{fundamental}.
\end{align}
From this equation we see that
\begin{align*}
\mathcal{A}L\mathcal{A}\sp{-1}&=A+B\tau,\label{char1}
\end{align*}
and consequently
\[
  \Tr L\sp{n}=\Tr(A+B\tau)\sp{n}.
\]
Thus we may relate the character of the automorphism $\sigma$ to
the trace of $A+B\tau$. For curves with many
symmetries we may evaluate $\tau$ via character theory. Streit
\cite{s01} examines period matrices and representation
theory. From (\ref{fundamental}) we obtain the \textbf{algebraic
  Riccati equation}
\begin{equation*}
\tau B\tau+\tau A-D\tau-C=0.
\end{equation*}
This fixed point equation may be understood as a
manifestation of the Riemann bilinear relations. The strategy then
is to calculate for appropriate automorphisms the matrices $M$ and
$L$ with the aim of restricting the form of the period matrix.
Such an approach goes back to Bolza \cite{b88} who studied the
consequences of a surface admitting many automorphisms in the
$g=2$ setting while more recently Schiller \cite{s68,s69} and
Schindler \cite{s93} have looked at the period matrices of
hyperelliptic curves using this approach. When a curve has a real
structure the period matrix is also constrained and Buser and
Silhol \cite{bs01} construct period matrices for various real
hyperelliptic curves. We also note that Bring's curve, the genus
$4$ curve with maximal automorphism group $S_5$, has a period
matrix essentially determined by its symmetries \cite{rr92}.

\section{Klein's curve} Klein's quartic curve may be expressed as
the plane algebraic curve
\begin{displaymath}
  x^3 y + y^3 z + z^3 x = 0.
\end{displaymath}
It is of genus three and achieves the Hurwitz bound of 168 on the
order of its (conformal) automorphism group $G$, which here is
isomorphic to $PSL(2,7)$. Klein's curve also has a real structure
possessing the antiholomorphic involution $[x,y,z]\mapsto
[\bar{x},\bar{y},\bar{z}]$. The entire symmetry group has 336
elements and is $PGL(2,7)$.

The affine projection ($(x,y)\leftrightarrow[x,y,1]$)
\begin{displaymath}
  x^3 y + y^3 + x = 0
\end{displaymath}
has 9 branch-points in the $x$-plane at 0, $\infty$ and the points
of a regular septagon centred on 0. The birational transformation
\begin{align*}
  (x,y) \mapsto (t,s) &= \left(
    1+\frac{x^3}{y^2},
    -\frac{x}{y}
  \right), \qquad
  (t,s) \mapsto (x,y) = \left(
    \frac{t-1}{s^2},
    \frac{1-t}{s^3}
    \right),
\end{align*}
yields the curve $$
  s^7 = t(t-1)^2
$$ which has 3 branch points on the $t$ plane located at 0, 1 and
$\infty$. In order to make the order three generator of $G$ more
manifest we perform a fractional-linear transformation sending these
branch-points to the cube roots of unity. Setting
$\rho=\exp(2\pi\mi/3)$ the transformation
\begin{align*}
  (t,s)\mapsto (z,w) = \left(
    \frac{t+\rho^2}{\rho t + \rho^2},
    \frac{s(\rho^2-1)}{\rho t + \rho^2}
  \right), \qquad
  (z,w)\mapsto (t,s) = \left(
    \frac{\rho^2(z-1)}{1-\rho z},
    \frac{w}{1-\rho z}
  \right)
\end{align*}
then leads to the curve
\begin{equation}
  w^7 = (z-1)(z-\rho)^2(z-\rho^2)^4.
\end{equation}
We will denote by $\phi:(x,y)\rightarrow(z,w)$ the composition of
these maps.

We choose a basis of differentials $ \{\omega_1,\omega_2,\omega_3\}$
where
\begin{align*}
  \label{eq:diffs-xy}
  \omega_1&=\frac{x\dif x}{x^3+3y^2},\qquad
  (\phi^{-1})^{*} \omega_1
    = \frac{\rho-1}{7}(z-\rho)(z-\rho^2)^2\frac{\dif z}{w^5},\\
  \omega_2&=\frac{y\dif x}{x^3+3y^2},\qquad
  (\phi^{-1})^{*} \omega_2
    = \frac{-1-2\rho}{7}(z-\rho)(z-\rho^2)^3\frac{\dif z}{w^6}, \\
  \omega_3&=\frac{\dif x}{x^3+3y^2},\qquad
  (\phi^{-1})^{*} \omega_3
    = \frac{2+\rho}{7}\frac{\dif z}{w^3}.
\end{align*}

\subsection{Symmetries} For our purposes we will focus on elements of
the symmetry group of orders 2, 3 and 7 as well as the
antiholomorphic involution. We describe these in the order that we
will employ them when determining the period matrix.

\subsubsection{Order 3 cyclic automorphism}
The projective representation of Klein's curve has the manifest
cyclic symmetry $[x,y,z] \mapsto [y,z,x]$ which yields the order
three automorphism
\[
  (x,y) \mapsto
    \left(\frac{y}{x},\frac{1}{x}\right), \qquad
  (z,w) \mapsto
    \left(\rho^2 z,
      \frac{-\rho^2(z-1)(z-\rho)(z-\rho^2)^2}{w^3}\right).
\]
It gives the simple cyclic rotation on (the pullback of) holomorphic
differentials
\begin{align*}
  \omega_1  \mapsto \omega_2, \ \
  \omega_2  \mapsto \omega_3, \ \
  \omega_3  \mapsto \omega_1.
\end{align*}

\subsubsection{Antiholomorphic involution} This leads to
\begin{displaymath}
  (z,w) \mapsto \left(\frac{1}{\bar{z}}, -\rho\frac{\bar{w}}{\bar{z}}\right).
\end{displaymath}
The symmetry simply acts as complex conjugation on the
differentials.

\subsubsection{Order 7 automorphism} Set $\zeta = \exp(2\pi\mi/7)$.
Klein's curve has an order 7 automorphism
\[
  (x,y) \mapsto (\zeta^5 x, \zeta^4 y),\qquad
  (z,w) \mapsto (z, \zeta w).
\]
The pullback action on differentials is
\begin{align*}
  \omega_1 \mapsto \zeta^2 \omega_1, \ \
  \omega_2 \mapsto \zeta^1 \omega_2, \ \
  \omega_3 \mapsto \zeta^4 \omega_3.
\end{align*}

\subsubsection{Holomorphic involution} This symmetry is really the
square of an order 4 automorphism in $G$ but has a simpler
description. As it suffices to fully restrict the form of the
period matrix we focus on this. The idea is that if we look at the
set of real solutions in the correct coordinates it has an obvious
rotational symmetry of order 2; this extends to the space of
complex solutions and gives the involution.

\begin{figure}
  \centering
  \includegraphics{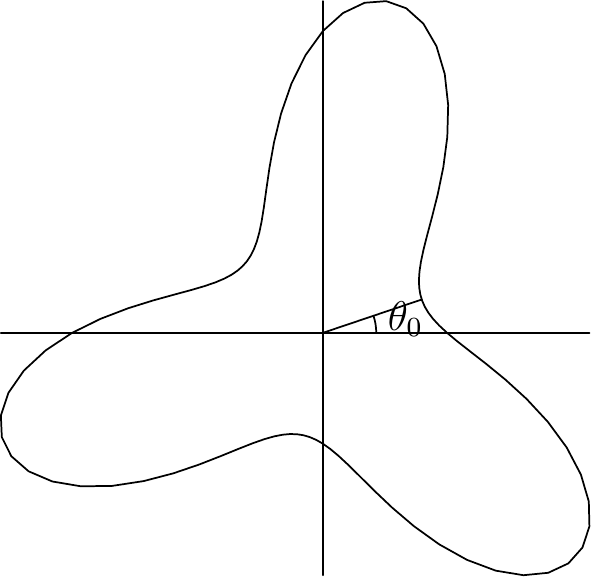}
  \caption{Section of real cone at $\tilde{z}=1$ after rotation.}
  \label{fig:z1-plot}
\end{figure}

So we proceed in stages from the original $x^3y+y^3z+z^3x=0$ form.
The set of real solutions of this homogeneous equation forms a
cone. If we rotate axes so that the axis of this cone, the
$(1,1,1)$ direction, becomes the $\tilde z$ axis of cylindrical
polar coordinates the curve takes the form
\begin{displaymath}
  4{\tilde z}^4 + 6r^2{\tilde z}^2 - 3r^4 - 2\sqrt{14} r^3 {\tilde z}\cos[3(\theta-\theta_0)]=0,
\end{displaymath}
where
$
  \theta_0 = \frac{\pi}{4} - \frac{1}{3}\tan^{-1} (3\sqrt{3})
$. This has an obvious threefold symmetry on rotating about
${\tilde z}$. See Figure \ref{fig:z1-plot} for a section of the
rays.
The involution is a rotation of $\pi$ radians about the
origin. Specifically
\begin{displaymath}
  (r,\theta,\tilde{z}) \mapsto (r,\theta+\pi,-\tilde{z}).
\end{displaymath}
Returning to the original coordinates,
\begin{equation}
  \begin{pmatrix}
    x\\y\\z
  \end{pmatrix}
  \mapsto
  \begin{pmatrix}
    \alpha & \beta & \gamma \\
    \beta & \gamma & \alpha \\
    \gamma & \alpha & \beta
  \end{pmatrix}
  \begin{pmatrix}
    x\\y\\z
  \end{pmatrix}
\end{equation}
where
\begin{align}
  3\alpha &= \cos\left(\frac{2}{3}\tan^{-1}(3\sqrt{3})\right)
      -\sqrt{3}\sin\left(\frac{2}{3}\tan^{-1}(3\sqrt{3})\right)-1,
      \nonumber\\
  3\beta &= -2\cos\left(\frac{2}{3}\tan^{-1}(3\sqrt{3})\right)-1,
    \label{eq:invol-algs}\\
  3\gamma &= \cos\left(\frac{2}{3}\tan^{-1}(3\sqrt{3})\right)
      +\sqrt{3}\sin\left(\frac{2}{3}\tan^{-1}(3\sqrt{3})\right)-1. \nonumber
\end{align}
In addition, a calculation shows that
\[
  \omega_1 \mapsto \alpha\omega_1 + \beta\omega_2 + \gamma\omega_3, \qquad
  \omega_2 \mapsto \beta\omega_1 + \gamma\omega_2 + \alpha\omega_3, \qquad
  \omega_3 \mapsto \gamma\omega_1 + \alpha\omega_2 + \beta\omega_3.
\]
For future reference we note that the $\alpha, \beta, \gamma$ are
not independent but satisfy
\begin{equation}\label{eq:abc}
  \alpha\gamma = \beta(\beta+1),\qquad
   \beta^2 = (\alpha+1)(\gamma+1).
\end{equation}

\section{A Canonical homology basis}

\begin{figure}
  \centering
  \includegraphics{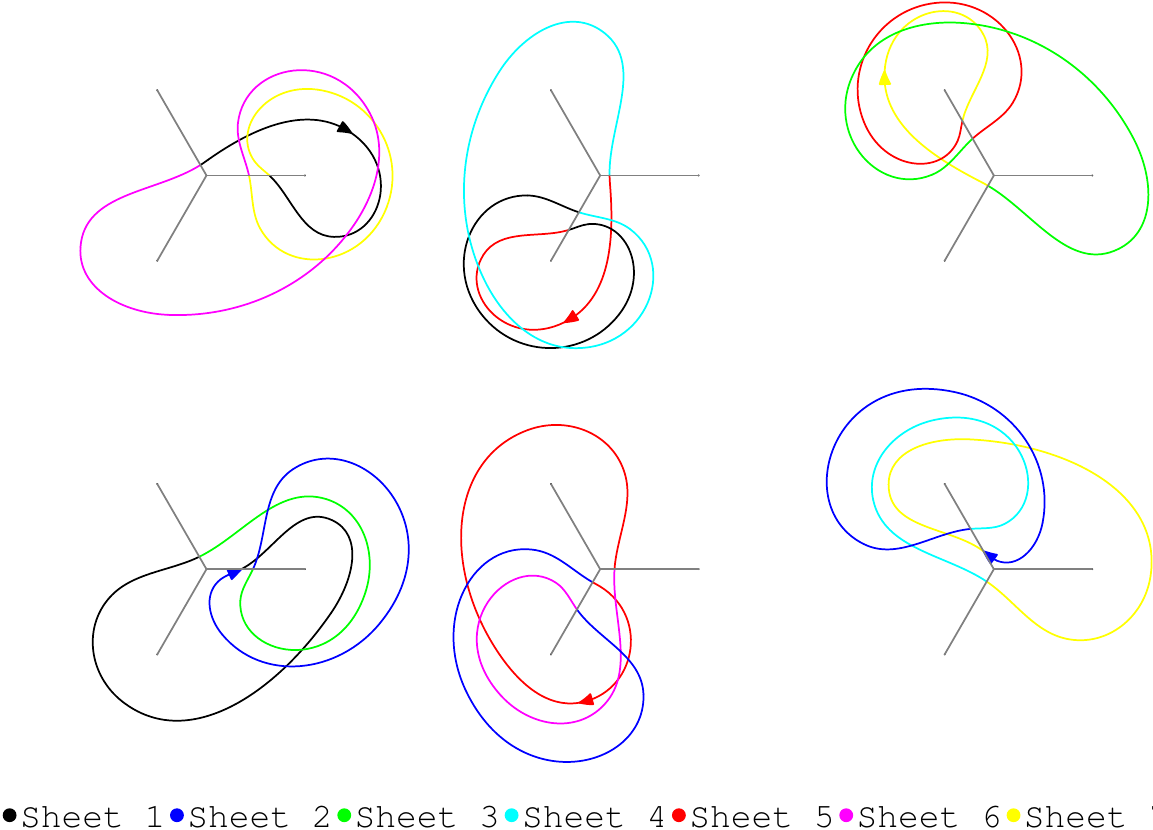}
  \caption{Homology basis in $(z,w)$ coordinates}
  \label{fig:homol-zw}
\end{figure}

We wish to construct a homology basis so that its transformation
under most of these symmetries is as simple as possible. The curve
has genus 3, and so we are looking for three $\mathfrak{a}$-cycles
and $\mathfrak{b}$-cycles. Thus a natural choice is to specify
$\mathfrak{a}_1$ and $\mathfrak{b}_1$ and then define the others
as images of the order 3 symmetry. This is most easily
accomplished in the $(z,w)$ space where there are only three
branch-points to choose from  and the automorphism is a simple
rotation in the $z$-plane. We need to know the monodromy around
the three branch points. Labelling the sheets by the $w$ value at
$z=0$ (i.e. for $k=0,\dots,6$ sheet $k$ is $(z,w) =
(0,\exp[{\mi\pi}(6k-1)/{21}])$) we find the monodromy about this
base point is (for each branch point) a simple constant shift of
sheets. The exact change is given in Table \ref{tab:monodromy}.
\begin{table}
  \centering
  \begin{tabular}{c|c}
    Branch point $z$ & Effect on sheet $k$ \\
    \hline
    $0$ & $k\mapsto k+1\pmod{7}$ \\
    $\rho$ & $k\mapsto k+2 \pmod{7}$ \\
    $\rho^2$ & $k\mapsto k+4 \pmod{7}$
  \end{tabular}
  \caption{Monodromy from base point 0 around each branch in a positive direction}
  \label{tab:monodromy}
\end{table}
Using this data we find that upon defining $\mathfrak{a}_1$ to
start at sheet 0 and proceed clockwise 3 times around $z=1$ and
then 1 time around $z=\rho^2$ we will obtain two additional
nonintersecting cycles upon applying the order 3 automorphism.
Further, by a simple shifting of the sheet, if $\mathfrak{b}_1$
starts on sheet 2 and the others are derived using the order 3
symmetry we obtain a canonical basis. See Figure
\ref{fig:homol-zw} for illustration (this is based on
\CyclePainter output).
It is now straightforward to transform this homology basis to the
other coordinate systems as needed. This homology basis will be
used to prove Theorem 1.

\subsection{Action of symmetries on the homology basis}
Having defined our homology basis and the symmetries it is a simple
matter to apply one to the other and write down the action of each
automorphism on the homology. It is here we have implemented code to
do this. With $\gamma =\left(
    \mathfrak{a}_1 , \mathfrak{a}_2 , \mathfrak{a}_3 ,
    \mathfrak{b}_1 , \mathfrak{b}_2 , \mathfrak{b}_3
  \right)^T$ we calculate intersection numbers to determine the matrix
$M$ introduced earlier. If we are considering the automorphism
$\sigma$ with
$
  \sigma_* : \gamma_i \mapsto M_{ij} \gamma_j,
$
and $\left<\cdot,\cdot\right>$ represents the intersection number then
\begin{displaymath}
\left<\sigma_*(\gamma_i), \gamma_k\right>
  = M_{ij} \left<\gamma_j, \gamma_k\right>
  = M_{ij} J_{jk},
\end{displaymath}
so $M_{ij} = \left<\sigma_*(\gamma_i), \gamma_k\right> J_{kj}^{-1}$
which is a trivial combination of intersection numbers. The results
are
\begin{itemize}
\item Order 3 automorphism. No computation is needed here. By
  definition of the homology this cyclically permutes the $A_i$ and
  $B_i$ independently.
  \begin{displaymath}
    M =
    \begin{pmatrix}
      0 & 1 & 0  &  0 & 0 & 0 \\
      0 & 0 & 1  &  0 & 0 & 0 \\
      1 & 0 & 0  &  0 & 0 & 0 \\
      0 & 0 & 0  &  0 & 1 & 0 \\
      0 & 0 & 0  &  0 & 0 & 1 \\
      0 & 0 & 0  &  1 & 0 & 0 \\
    \end{pmatrix}.
  \end{displaymath}

\item Antiholomorphic involution.
  \begin{displaymath}
    M =
    \begin{pmatrix}
      0 & 0 & 0   &  -1 & 0 & 0 \\
      0 & 0 & 0   &  0 & -1 & 0 \\
      0 & 0 & 0   &  0 & 0 & -1 \\
      -1 & 0 & 0  &  0 & 0 & 0 \\
      0 & -1 & 0  &  0 & 0 & 0 \\
      0 & 0 & -1  &  0 & 0 & 0 \\
    \end{pmatrix}.
  \end{displaymath}

\item Order 7 automorphism. Although this acts simply on paths
    by just shifting sheets a fairly complicated matrix
    results.
  \begin{displaymath}
    M =
    \begin{pmatrix}
      1 & 0 & -1 & 1 & 0 & -1 \\
      0 & 0 & 0 & 0 & 1 & 0 \\
      -1 & 0 & -1 & 0 & 1 & 0 \\
      -1 & 0 & 0 & 0 & 1 & 0 \\
      0 & -1 & -1 & 1 & 0 & 0 \\
      1 & 0 & 0 & 0 & 0 & -1
    \end{pmatrix}.
  \end{displaymath}
  However, since we know that $\mathfrak{b}_i$ \emph{is} simply $\mathfrak{a}_i$ shifted by
  some sheet number, we know that some power of this automorphism will
  map $\mathfrak{a}_i$ to $\mathfrak{b}_i$. Indeed $M$ itself takes
  $\mathfrak{a}_2$ to $\mathfrak{b}_2$;
  $M^2$ takes $\mathfrak{a}_1$ to $\mathfrak{b}_1$; and $M^4$ takes
   $\mathfrak{a}_3$ to $\mathfrak{b}_3$.
\item Order 2 involution. This has a remarkably simple effect on the
  homology, given its complexity as a transformation.
  \begin{displaymath}
    M =
    \begin{pmatrix}
      0 & 0 & -1  &  0 & 0 & 0 \\
      0 & -1 & 0  &  0 & 0 & 0 \\
      -1 & 0 & 0  &  0 & 0 & 0 \\
      0 & 0 & 0   & 0 & 0 & -1   \\
      0 & 0 & 0   & 0 & -1 & 0  \\
      0 & 0 & 0   & -1 & 0 & 0
    \end{pmatrix}.
  \end{displaymath}
\end{itemize}
We note in passing that these matrices were produced by applying
the \texttt{extcurves} intersection routines in the form
\begin{verbatim}
> find_homology_transform(curve, homSrc, homDst):
\end{verbatim}
for appropriate choices of the parameters, where
\texttt{homSrc}$=\{\mathfrak{a}_i,\mathfrak{b_i}\}$ and
\texttt{homDst}=$\{\sigma_\ast(\mathfrak{a}_i),\sigma_\ast(\mathfrak{b_i})\}$
and $\sigma$ is the automorphism being considered.

\section{The period matrix}
We now have the matrices $M$ and $L$ necessary to employ
(\ref{fundamental}). We shall briefly remark on the restrictions
that each symmetry results in.
\subsection{Order 3 symmetry}
This symmetry constrains the $\mathfrak{a}$-periods and
$\mathfrak{b}$-periods separately, but in an identical manner, so
we only comment on the former. Here the matrix relation is
\begin{displaymath}
  \begin{pmatrix}
    0 & 1 & 0 \\
    0 & 0 & 1 \\
    1 & 0 & 0
  \end{pmatrix}
  \mathcal{A}
  =
  \mathcal{A}
  \begin{pmatrix}
    0 & 0 & 1 \\
    1 & 0 & 0 \\
    0 & 1 & 0
  \end{pmatrix},
\end{displaymath}
which results in
\begin{displaymath}
  \mathcal{A} =
  \begin{pmatrix}
    X & Y & Z \\
    Y & Z & X \\
    Z & X & Y
  \end{pmatrix}.
\end{displaymath}

\subsection{Antiholomorphic involution}
This determines the $\mathfrak{b}$-periods in terms of the
$\mathfrak{a}$-periods. Calling the symmetry $\sigma$, we have
\begin{align*}
  \mathcal{B}_{ij} &= \int_{\mathfrak{b}_i} \omega_j
    = -\int_{\sigma(\mathfrak{a}_i)} \omega_j
    = -\int_{\mathfrak{a}_i} \sigma^*(\omega_j)
    = -\int_{\mathfrak{a}_i} \overline{\omega}_j
    = -\overline{\int_{\mathfrak{a}_i} \omega_j}
    = -\overline{\mathcal{A}}_{ij}.
\end{align*}

\subsection{Order 7 automorphism}

This symmetry sends $\mathfrak{a}_2$ to $\mathfrak{b}_2$ so it tells
us
\begin{align*}
  -\overline{Y}= \zeta^2 Y,\qquad
  -\overline{Z}= \zeta   Z,\qquad
  -\overline{X}= \zeta^4 X.
\end{align*}
This fixes the argument of all three numbers, up to sign. So if we
merely constrain $r_1, r_2, r_3$ to be real then we may write
\begin{align*}
  X &= r_1 \exp(-\pi\mi/14), &
  Y &= r_2 \exp(-11\pi\mi/14), &
  Z &= r_3 \exp(-9\pi\mi/14)
\end{align*}
(The choice has actually been made so that all $r_i$ are positive.)

\subsection{Holomorphic involution} Here we find that
\begin{equation}
  \label{eq:invol-eqs}
  -
  \begin{pmatrix}
    Z&X&Y\\ Y&Z&X\\ X&Y&Z
  \end{pmatrix}
  =
  \begin{pmatrix}
    X&Y&Z\\ Y&Z&X\\ Z&X&Y
  \end{pmatrix}
  \begin{pmatrix}
    \alpha&\beta&\gamma\\
    \beta&\gamma&\alpha\\
    \gamma&\alpha&\beta
  \end{pmatrix}.
\end{equation}
While this naively looks like three independent equations the
relations (\ref{eq:abc}) in fact mean that there is only one
independent (complex) equation in \eqref{eq:invol-eqs}. For
example, if we multiply the first equation $\alpha X + \beta Y +
\gamma Z = -Z$ in  \eqref{eq:invol-eqs} by $(\beta+1)/\alpha$ we
obtain (up to rewriting) the second equation
  $(\beta+1)X + \gamma Y + \alpha Z = 0$.
Now as at this stage we only have three real parameters $r_1, r_2,
r_3$ left this one complex equation is sufficient to determine
(say) $r_1$ and $r_2$ in terms of $r_3$ and hence the matrix of
periods up to an overall real multiple. Indeed, this real multiple
cancels when calculating the the Riemann form of the period
matrix, but as we will see in the next subsection even this scalar
parameter is explicitly calculable.

Take the first entry of \eqref{eq:invol-eqs} as our starting point;
\begin{displaymath}
  \alpha X + \beta Y + \gamma Z = -Z.
\end{displaymath}
With $r_1 = \mu r_3$, $r_2 = \nu r_3$ we find a quadratic solution for
$\mu, \nu$ in terms of the algebraic coefficients. Their minimal
polynomials over $\mathbb{Q}$ are simply
\[
  \mu^3 + \mu^2 - 2\mu - 1 = 0, \qquad
  \nu^3 - 2\nu^2 - \nu + 1 = 0.
\]
In fact each of these have three real roots, and split over
$\mathbb{Q}(\zeta)$. The solution is seen to be
\[
  \mu = \zeta + \zeta^{-1}, \qquad
  \nu = 1 + \zeta + \zeta^{-1}.
\]

\subsection{The final free parameter}

The final free parameter can be computed most readily from the
integral in the $(t,s)$ space. The cycle $\mathfrak{a}_1$ can be
deformed to a path traversing $t=0$ to $t=1$ and then back again on
a different sheet. Along this segment there is a sheet where $s$
remains real, and all other sheets are simply related by
multiplication with some power of $\zeta$. Thus for example
\begin{align*}
  Z &=\int_{\mathfrak{a}_1} \omega_3
    = \int_{\mathfrak{a}_1} \frac{\dif t}{7s^3}
    = \frac{\zeta^k-\zeta^l}{7}\int_0^1 \frac{\dif t}{s^3}
%    = \frac{\zeta^k}{7}\int_0^1 \frac{\dif t}{t^{3/7}(t-1)^{6/7}}
%       + \frac{\zeta^l}{7}\int_1^0 \frac{\dif
%       t}{t^{3/7}(t-1)^{6/7}}\\
    = \frac{\zeta^k-\zeta^l}{7}\int_0^1 \frac{\dif t}{t^{3/7} (t-1)^{6/7}}
    = \frac{B(\frac{4}{7}, \frac{1}{7})}{7}(\zeta^k -\zeta^l).
\end{align*}
All that remains is to determine the integers $k$ and $l$. This is a
trivial matter of analytically continuing the path $\mathfrak{a}_1$
over to somewhere on the segment $(0,1)$ and noting the argument of
the result. We discover
\begin{displaymath}
  Z = \frac{B(\frac{4}{7}, \frac{1}{7})}{7}(\zeta^{-1} - 1).
\end{displaymath}
Using this to express $r_3$ in terms of $\Gamma$-functions rather
than $B$-functions
\begin{displaymath}
  r_3 =
  2\sin(\pi/7)\frac{\Gamma(\frac{8}{7})\Gamma(\frac{4}{7})}
  {\Gamma(\frac{5}{7})}.
\end{displaymath}

\subsection{The Riemann period matrix}

The Riemann period matrix is
$
  \tau = \mathcal{B}\mathcal{A}^{-1}.
$ Bringing our results together yields (\ref{eq:kleinperiod}).

\section{The relation to the period matrix of Rauch and Lewittes}

\subsection{Hyperbolic model of Klein's curve}

Rauch and Lewittes \cite{rl69} have already described a canonical
homology basis for this curve and calculated the associated Riemann
period matrix. In general finding the symplectic transformation
between two period matrices is a hard problem. The best approach is
usually to work with the homology bases involved and relate those.

Rauch and Lewittes' description is based on the expression of Klein's
curve as a quotient of the hyperbolic disc, so in order to compare the
two results we will need some understanding of that model.

\begin{figure}
  \centering
  \includegraphics{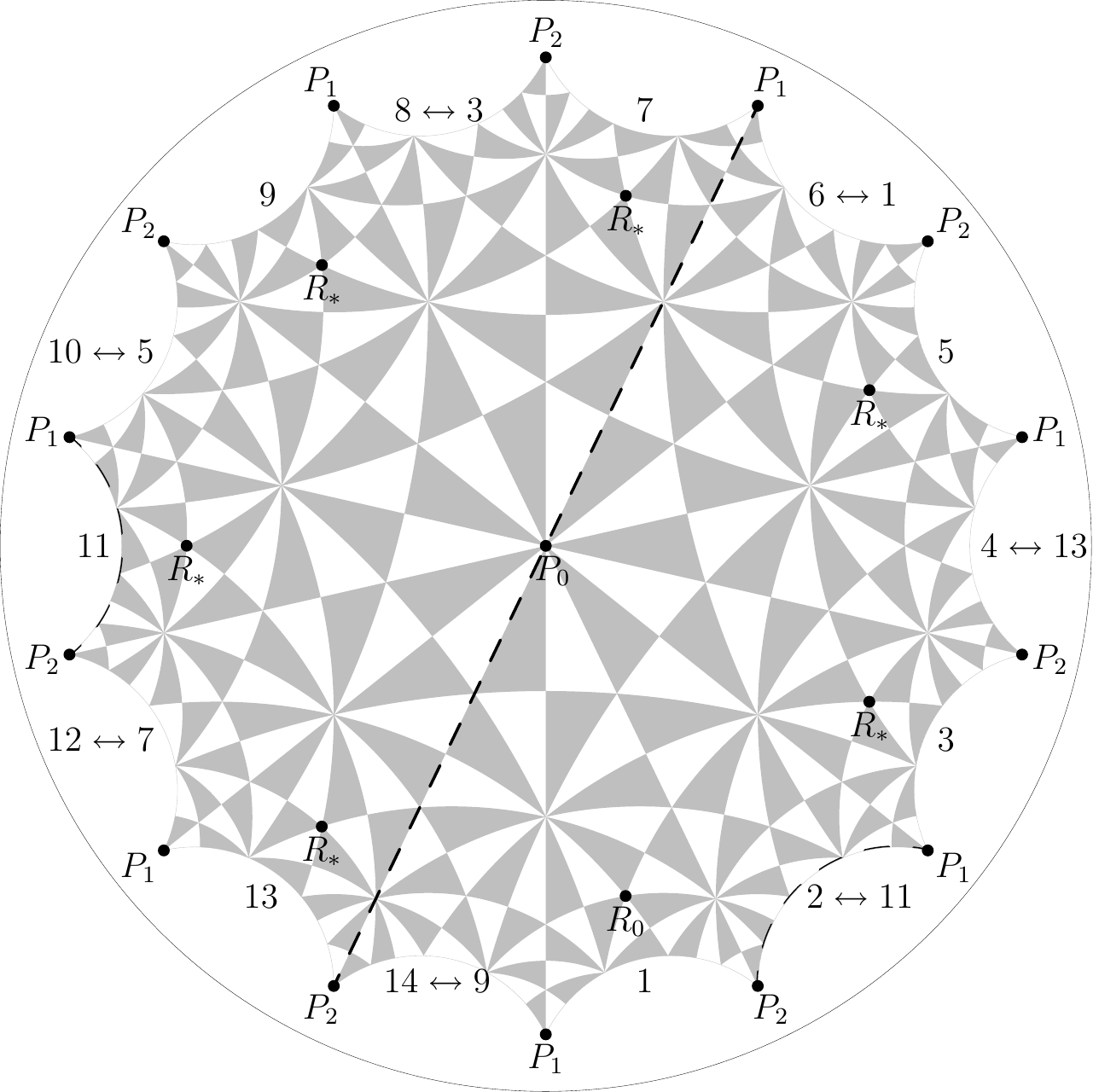}
  \caption{Hyperbolic disc model of Klein's curve.}
  \label{fig:hyperbolic}
\end{figure}

The hyperbolic model amounts to a regular 14-gon centred in the
disc, and the identification of vertices and pairs of sides around
this as indicated in Figure \ref{fig:hyperbolic}. It may be tiled
by 336 triangles with angles $(2\pi/7, \pi/3, \pi/2)$, each of
which can form a fundamental domain of the symmetry group. This
group is then manifest as rotations about any triangular vertex
together with reflections in any geodesic line consisting of
triangular edges.

Rauch and Lewittes' homology basis is described in terms of paths
back and forth between $P_1$ and $P_2$ along prescribed edges.
Taking positive numbers to indicate anticlockwise traversal (about
the centre) and negative the reverse, the basis is explicitly
\begin{align*}
  \mathfrak{a}'_1 &= 1 - 4 - 7 - 9, \qquad
  \mathfrak{a}'_2 = -4 - 9, \qquad
  \mathfrak{a}'_3 = -4 - 5, \\
  \mathfrak{b}'_1 &= 2+3+4+5, \qquad
  \mathfrak{b}'_2 = -3+7, \qquad
  \mathfrak{b}'_3 = 3-5.
\end{align*}
(So for example $\mathfrak{a}'_1$ consists of moving from $P_1$ to
$P_2$ along side 1, back to $P_1$ along 7, to $P_2$ again along 4 and
finally back to the start along 9. Although there is ambiguity in the
order of edges taken all the resulting paths are homologous.)

\subsection{Identification of the two models}

We now wish to identify the two models of Klein's curve so we can
express Rauch and Lewittes' homology basis in some coordinate
plane. Working in $(t,s)$ coordinates we are seeking two meromorphic
functions $t, s$ on the hyperbolic disc with the property that
\[
  s(p)^7 = t(p)(t(p)-1)^2
\]
everywhere. If we could describe such functions well enough then
we could write down the coordinate projection of the hyperbolic
basis. There is obviously no unique choice for such functions
since the push-forward of any pair under an automorphism would in
general give a different pair. However, we can exploit this to
produce some $t,s$ with useful properties.

Consider the subgroup $H$ of the full automorphism group generated by
\begin{align*}
  Z : (t,s) &\mapsto (t,\zeta s), &
  R : (t,s) &\mapsto \left(1-\frac{1}{t},\frac{t-1}{s^3}\right), &
  c : (t,s) &\mapsto (\bar{t}, \bar{s}).
\end{align*}
These satisfy the relations $Z^7 = R^3 = c^2 = 1$, $RZR^{-1} = Z^4$,
$cRc^{-1} = R$, $cZc^{-1} = Z^{-1}$ which allow us to write
\[
  H = \{ c^i Z^j R^k | i = 0,1;\ j = 0,\ldots,6;\ k=0,1,2 \},
\]
and hence $H$ has 42 elements.

It is easy to verify that the full symmetry group of Klein's curve
(which is isomorphic to $PGL(2,7)$) has eight subgroups of order
42, all of which are conjugate. So if we can find an appropriate
such subgroup on the hyperbolic model (one which interacts well
with Rauch and Lewittes' basis ideally) then by applying some
automorphism we may assume without loss of generality that this is
our group. This is helpful because, for example, the
transformation $Z$ has a simple description in terms of the $s$,
$t$ coordinates.

To obtain the hyperbolic group, we start with rotations of order 7
about $P_0$ and adjoin the rotations of order three about the
indicated points $R_{\bullet}$. This in itself is a group, and has
order 21 (1 identity + 6 rotations of order 7 + $7\times 2$
rotations of order 3). Let $\hat{H}$ be the result of adjoining a
reflection in a diameter through $P_1$ to this group. It has the
required 42 elements and can be made the counterpart of $H$ in the
hyperbolic model.

Specifically, suppose we found an arbitrary pair of functions $t,s$
with the required algebraic relations. This induces a hyperbolic
subgroup corresponding to $H$. If this is not $\hat{H}$ then they are
at least conjugate so $\hat{H} = g^{-1} H g$. But then if we instead
use the meromorphic pair $t\circ g, s\circ g$ then $H$ will induce
$\hat{H}$ itself. In fact, by similar means, remaining within the group
$\hat{H}$ we can demand more.

We start with $c$. This fixes an entire line pointwise in $H$. The
only elements of $\hat{H}$ with this property are reflections in
the diametric geodesics, so $c$ must correspond to one of those.
But we can do better: conjugating with some central rotation in
the same way we did before we can guarantee that $c$ corresponds
to the reflection that fixes edge 2 as indicated in Figure
\ref{fig:hyperbolic}. Consequently we have guaranteed that the
functions $s(p)$, $t(p)$ are real on edge 2.

Having made this demand, we discover that the possible images for $R$
in $\hat{H}$ are severely limited. In order to satisfy $cRc^{-1} = R$,
$R$ must be some rotation about $R_0$ in $\hat{H}$. Finally, in order
to satisfy $RZR^{-1} = Z^4$ we actually need $R$ to correspond to the
anticlockwise rotation of $2\pi/3$ about $R_0$.

Now, $Z$ must correspond to one of the central rotations. So the
fixed points of $Z$, namely $t=0,1,\infty$ must correspond to
$P_0, P_1, P_2$ in some order. But since rotations about $R_0$
permute these we can conjugate again to demand that $t(P_0) =
\infty$. It is easily checked that this disrupts neither the group
$\hat{H}$ nor our previous choice of $c$ (essentially because $R c
R^{-1} = c^{-1} = c$). Before abandoning this issue we note that
$R(\infty,\bullet) = (1,\bullet)$ and that the corresponding
element of $\hat{H}$ sends $P_0$ to $P_1$ so $t(P_1) = 1$, $t(P_2)
= 0$.

At this stage we know that traversing side 2 (or 11) is equivalent
to travelling from $t=0$ to $t=1$ along the real axis. If we knew
that $Z$ ($s\mapsto \zeta s$) corresponded to a central rotation by
$2\pi j/7$ then we could deduce the phase of $s$ along any numbered
edge in terms of $j$. Specifically along edge $2k$
\begin{displaymath}
  s \in \zeta^{(k-1)j^{-1}} \mathbb{R},
\end{displaymath}
where the inverse $j^{-1}$ is taken (mod 7). The odd edges are obtained
by their identification with even edges. This would be enough to
completely determine any path expressed, as Rauch and Lewittes do,
by which numbered edges should be traversed in the hyperbolic
model.

So our final task is to identify $j$ or equivalently which central
hyperbolic rotation our coordinate transformation $Z$ corresponds
to. The solution is provided by the hyperbolic structure near the
point $P_1$. The idea is that $P_1$ is the centre of a septagon and if
you go towards it on some edge and away on another then the angle
between these paths corresponds to how many branch cuts you would
cross doing the same thing in coordinate space.  Requiring the same
phase for $s$ from both processes is enough to fix $j$.

\begin{figure}
  \centering
  \includegraphics{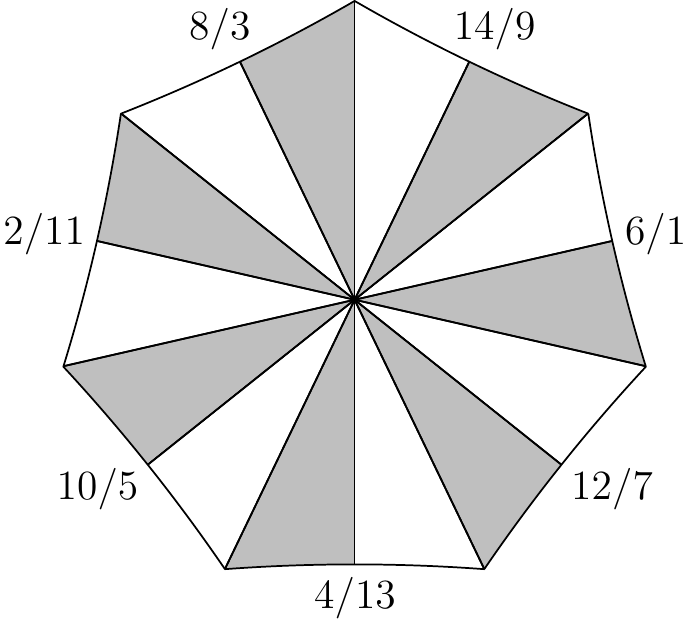}
  \caption{How edges come together near $P_1$}
  \label{fig:branches}
\end{figure}

More precisely, we can deduce from Figure \ref{fig:hyperbolic} how
the numbered edges are laid out around $P_1$ and $P_2$. For example
the bottom vertex shows us that an anticlockwise rotation of
$2\pi/7$ about $P_1$ will take us from edge 1 to edge 14. Putting
all this information together, we obtain Figure \ref{fig:branches}.

Now notice the unlabelled ``spokes'' in the diagram. We can think
of these as branch cuts in coordinate space. If a branch cut was
some ray from $t=1$ on all sheets, then in hyperbolic space these
would be represented by a single line and its images under the
sheet-changing operator $Z$ (or some power). These lines won't
necessarily be geodesic, but they \emph{will} have a well-defined
direction emerging from $P_1$ and be related by some $2\pi/7$
rotation in Figure \ref{fig:branches}. Thus (for example by
considering a small neighbourhood around $P_1$) we may as well
consider the unlabelled ``spokes'' as the branch cuts.

Now we can put these two pictures together. Suppose we start with both
$t$ and $s$ real. If we go around the branch point $t=1$ once
anticlockwise we discover that $s \in\zeta^2\mathbb{R}$. In doing so
we've crossed just one branch cut. In the hyperbolic picture what
we've done is start on edge 2 and go anticlockwise crossing one
spoke. So we are on edge 10. But this corresponds to $s \in
\zeta^{4j^{-1}} \mathbb{R}$. For these two statements to be consistent
$j=2$ and $Z$ corresponds to a central rotation of $4\pi/7$ and the
argument of $s$ along each edge is as in Table \ref{tab:args}.
\begin{table}
  \centering
  \begin{tabular}{l|ccccccc}
    Edge  & 2/11 & 4/13      & 6/1     & 8/3       & 10/5      & 12/7     & 14/9 \\
\hline
    Phase & 1    & $\zeta^4$ & $\zeta$ & $\zeta^5$ & $\zeta^2$ & $\zeta^6$ & $\zeta^3$
  \end{tabular}
  \caption{$s$ phase on each numbered edge}
  \label{tab:args}
\end{table}

Note that in the above there was an ambiguity over the direction of
paths -- an implicit assumption that anticlockwise was the same in
both models. This corresponds to a choice of orientation. Group
theoretically if we conjugate by $c$ it takes us between these two
choices (because $cZc^{-1} = Z^{-1}$). The choice made gives us a
symplectic transformation between the two homology bases, rather
than antisymplectic.

\subsection{Rauch and Lewittes' Homology in Coordinates}

\begin{figure}
  \centering
  \includegraphics{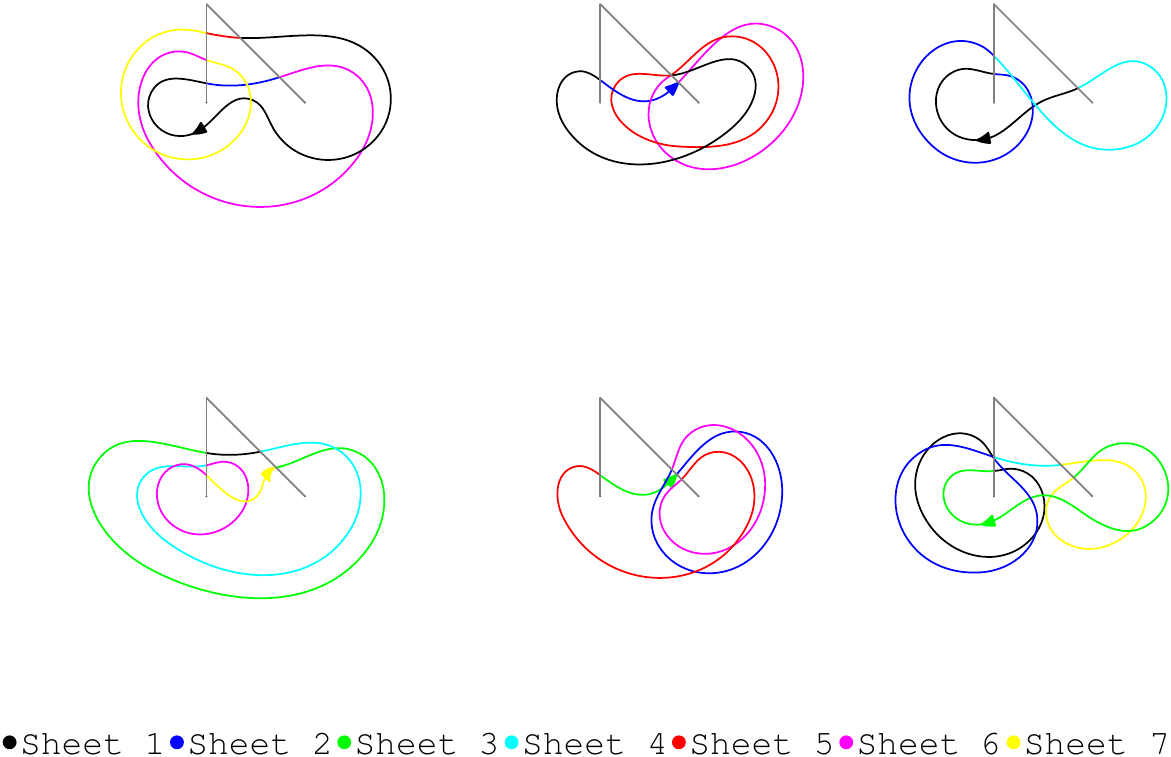}
  \caption{Rauch and Lewitte's homology basis in $(t,s)$ coordinates.}
  \label{fig:rauch-hom}
\end{figure}

With the above identification, we can simply read off which sheets
Rauch and Lewitte's homology basis uses in its repeated journeys
between $t=0$ and $t=1$. From that information we put together paths
that make the same journeys and so draw the basis. We obtain Figure
\ref{fig:rauch-hom}.

We can algorithmically check that these paths form a canonical
homology basis, and indeed we discover the correct intersection
matrix.

It is now a simple matter to calculate the symplectic
transformation moving from our basis to this one and the result is
the previously stated (\ref{eq:rlsymp}). The symplectic
equivalence between period matrices gives a nontrivial test of
this identification.

\section{Other period matrices}

Tretkoff and Tretkoff \cite{tretkoff}, Schindler \cite{s91},
Rodr\'iguez and Gonz\'alez-Aguilera \cite{rg95}, Yoshida
\cite{y99} and Tadokoro \cite{t08} have all described homology
bases for Klein's curve and calculated the resulting period
matrices. Here we relate their works with ours, beginning with
those presenting Klein's curve as a plane projective model.

\subsection{Yoshida's period matrix}

The symplectic transformation
\[
\begin{pmatrix}
  0 & 0 & 0 & -1 & 0  & 0  \\
  0 & 1 & 0 & 0  & 0  & -1 \\
  0 & 0 & 1 & 0  & -1 & 1  \\
  1 & 1 & 0 & 0  & 0  & -1 \\
  0 & 1 & 1 & -1 & 0  & 0  \\
  0 & 1 & 0 & 0  & 0  & 0
\end{pmatrix}
\]
takes our period matrix to that given by Yoshida,
\[
\frac{1}{14}
\begin{pmatrix}
  6\mi\sqrt{7} & 7 + 3\mi\sqrt{7} & 2\mi\sqrt{7} \\
  7 + 3\mi\sqrt{7} & 7 + 5\mi\sqrt{7} & 7 + \mi\sqrt{7} \\
  2\mi\sqrt{7} & 7 + \mi\sqrt{7} & 7+3\mi\sqrt{7}
\end{pmatrix}.
\]

\subsection{Tadokoro's period matrix}

Similarly the symplectic transformation
\[
\begin{pmatrix}
  0  & 0  & 0 & 1  & 0  & 0 \\
  0  & -1 & 1 & 0  & -1 & 1 \\
  0  & 1  & 0 & -1 & 0  & -1 \\
  -1 & 1  & 0 & -1 & 1  & 0 \\
  0  & 1  & 0 & -1 & 0  & 0 \\
  0  & 1  & 1 & -1 & 0  & 0
\end{pmatrix}.
\]
takes our period matrix to the period matrix given by Tadokoro,
\[
\frac{1}{8}
\begin{pmatrix}
  -1+3\mi\sqrt{7} & -3+\mi\sqrt{7} & 2+2\mi\sqrt{7} \\
  -3+\mi\sqrt{7} & \-1+3\mi\sqrt{7} & 2+2\mi\sqrt{7} \\
  2+2\mi\sqrt{7} & 2+2\mi\sqrt{7} & 4+4\mi\sqrt{7}
\end{pmatrix}.
\]

\subsection{Tretkoff and Tretkoff's period matrix}
In \cite{tretkoff} these authors develop the algorithm (named
after them) that determines an homology basis of a Riemann surface
and Klein's curve is one of the examples given there (p495-498).
This is the algorithm implemented in Maple's \texttt{algcurves}.
Here the author's reproduce the work of Hurwitz, Klein and Fricke,
and Baker with the choice of cycles given by those earlier works.
Finally they then relate these cycles to a canonical homology
basis and produce a period matrix. The symplectic transformation
\[
\begin{pmatrix}
0&0&-1&1&0&0\\ 0&1&1&
-1&0&0\\ 0&0&-1&0&0&0\\ -1&0&-1&1&
1&0\\ 0&0&0&0&1&0\\ 1&0&0&0&0&-1
\end{pmatrix}.
\]
takes our period matrix to
\[
\frac18
\begin{pmatrix}
7+3\mi\sqrt {7}&-3+\mi\sqrt {7}&-5-\mi\sqrt {7}\\
-3+\mi\sqrt {7}&-1+3\mi\sqrt {7}&-3+\mi\sqrt {7}\\
-5-\mi\sqrt {7}&-3+\mi\sqrt {7}&-1+3\mi\sqrt {7}
\end{pmatrix}.
\]
This is not quite the period matrix Tretkoff and Tretkoff give.
One finds that although the matrix given on the left-hand-side at
the top of p498 is correct, in the final evaluation an error has
been made in their $2\tau\sp3$ terms; correcting this yields the
above.

\subsection{Schindler's period matrix}
The period matrix of Klein's curve was also calculated in the
thesis \cite{s91}. Comparing the homology bases of our works
yields the following symplectic transformation
\begin{equation*} \label{eq:rlsympschindler}
  \begin{pmatrix}
   1&0&0&0&-1&0\\
   0&-1&-1&1&0&1\\
   0&-2&-1&2&0&1\\
   -1&0&-1&1&1&1\\
   2&1&1&-1&-1&-3\\
   -1&0&0&0&0&1
  \end{pmatrix}.
\end{equation*}
This leads to the period matrix
\begin{equation} \label{eq:tausympschindler}%(C+D\tau)(A+B\tau)^{-1}=
\frac1{14}\begin{pmatrix}
6\mi\sqrt {7}-14&21-5\mi\sqrt {7}&-7+3\mi
\sqrt {7}\\ 21-5\mi\sqrt {7}&-42+10\mi\sqrt {7}&14-6
\mi\sqrt {7}\\ -7+3\mi\sqrt {7}&14-6\mi\sqrt {7}&-7+
5\mi\sqrt {7}
\end{pmatrix}
%=\Pi_1'+\Pi_2\,i\sqrt{7}
.\end{equation} Here we find that the imaginary parts of this and
Schindler's period matrix agree but the real parts are in
disagreement. We believe a typographic error before the final
substitution led to an incorrect result.\footnote{We thank H.
Lange for correspondence on this point.} The result
(\ref{eq:tausympschindler}) should replace that given in
\cite[Exercise 11.15]{bl04} (the imaginary parts of these agree
however).

\subsection{The Rodr\'iguez and Gonz\'alez-Aguilera period matrix}
These authors consider the  classically studied KFT pencil of
curves
$$x^4+y^4+z^4+t(x^2y^2+y^2z^2+z^2x^2)=0$$
in terms of a hyperbolic model. For $t=-3(1\pm\mi\sqrt{7})/2$ this
is Klein's curve, for $t=0$ this is Fermat's quartic curve, and
for $t=-2$ the tetrahedron, and so the nomenclature. There are
various singular fibres in this pencil. Away from these the period
matrix takes the form
\[
\tau(t)
\begin{pmatrix}
  3&-1&-1 \\ -1&3&-1 \\ -1&-1&3
\end{pmatrix}
\]
and the Jacobian is isogenous to the product of three elliptic
curves. Klein's curve corresponds to $\tau(t)=(5+\mi\sqrt{7})/2$
and the symplectic transformation
\[
\begin{pmatrix}
  0&0&1&-1&0&0\\ 0&1&1&
-1&0&0\\ 0&0&1&0&0&0\\ 1&-1&1&-1&-
1&0\\ 0&2&0&-1&1&0\\ -1&-1&1&1&0&1
\end{pmatrix}.
\]
takes our period matrix to this.

\bibliographystyle{amsalpha}

\begin{thebibliography}{D{i}e91}

\bibitem[BBM]{bbm83}M, B. Babich,  A. I. Bobenko and V. B. Matveev,
\emph{Reductions of Riemann theta functions of genus $g$ to theta
functions of lesser genus, and symmetries of algebraic curves},
Dokl. Akad. Nauk SSSR  272  (1983),  no. 1, 13--17.

\bibitem[Ba]{b07}H.F. Baker, \emph{Multiply Periodic Functions},
Cambridge University Press, Cambridge, 1907.

\bibitem[BL]{bl04}Christina Birkenhake and Herbert Lange,
\emph{Complex Abelian Varieties} Second edition, Springer-Verlag,
Berlin 2004.

\bibitem[Bo]{b88}Oskar Bolza, \emph{
On Binary Sextics with Linear Transformations into Themselves},
Amer. J. Math. 10 (1887), no. 1, 47--70.

\bibitem[BE07]{bren07}
H.~W. Braden and V.~Z. Enolski, \emph{Monopoles, {C}urves and
{R}amanujan}, %Reported at Riemann {S}urfaces,
  %{A}nalytical and {N}umerical {M}ethods, {M}ax {P}lanck {I}nstititute
  %(Leipzig), 2007. Submitted,
To be published. arXiv: math-ph/0704.3939.

\bibitem[BE09]{bren09}
\bysame, \emph{On the tetrahedrally symmetric monopole}, arXiv:
\texttt{math-ph/0908.3449}

\bibitem[B]{b10} H.~W. Braden,  \emph{Cyclic Monopoles, Affine Toda
and Spectral Curves}, arXiv:1002.1216

\bibitem[BK]{BK}
{Integrability: the Seiberg-Witten and Whitham equations}, H. W.
Braden and I. M. Krichever (Eds.), Gordon and Breach Science
Publishers, Amsterdam, 2000. ISBN: 90-5699-281-3.

\bibitem[BS]{bs01}Peter Buser and Robert Silhol, \emph{Geodesics,
Periods and Equations of Real Hyperelliptic Curves}, Duke Math. J.
108 (2001) 211--244.

\bibitem[D]{d10}A. D'Avanzo,
\emph{On Charge 3 Cyclic monopoles}, Edinburgh Ph.D. Thesis, 2010.

\bibitem[FK]{fk}
Robert Fricke and  Felix Klein, \emph{Vorlesungen über die Theorie
der automorphen Funktionen. Band II: Die funktionentheoretischen
Ausführungen und die Andwendungen} Teubner, Stuttgart 1912.

\bibitem[G]{jg07} Jane Gilman, \emph{Canonical symplectic representations
 for prime order conjugacy classes of the mapping-class group},
  J. Algebra  318  (2007) 430--455.

\bibitem[HMM95]{hmm95}
N. ~J. Hitchin, N. ~S. Manton and  M. ~K. Murray, \emph{Symmetric
monopoles}, Nonlinearity \textbf{8} (1995), 661--692.

\bibitem[H]{hur} Adolf Hurwitz, \emph{Ueber einige besondere homogene lineare
Differentialgleichungen}, Math. Ann. 26 (1886), no. 1, 117--126.

\bibitem[L]{l99} Silvio Levy (Editor), \emph{The Eightfold Way: The
 Beauty of Klein's Quartic Curve}, Math. Sci. Research Inst. Publications.

\bibitem[RL]{rl69}Harry E. Rauch and J. Lewittes,
\emph{The Riemann surface of Klein with 168 automorphisms}, Problems
in analysis (papers dedicated to Salomon Bochner, 1969), pp.
297--308. Princeton Univ. Press, Princeton, N.J., 1970.

\bibitem[RR]{rr92}Gonzalo Riera and  Rub\'i E Rodriguez,
\emph{The period matrix of Bring's curve},  Pacific J. Math.  154
(1992),  no. 1, 179--200.

\bibitem[RG]{rg95}Rub\'i E Rodr\'iguez and V\'ictor Gonz\'alez-Aguilera,
\emph{Fermat's quartic curve, Klein's curve and the tetrahedron},
in Extremal Riemann surfaces (San Francisco, CA, 1995), 43--62,
Contemp. Math., 201, Amer. Math. Soc., Providence, RI, 1997.

\bibitem[JS-a]{s68}John Schiller, \emph{Riemann matrices for hyperelliptic
surfaces with involutions other than the interchange of sheets},
Michigan Math. J. 15 (1968) 283--287.

\bibitem[JS-b]{s69}John Schiller, \emph{Moduli for special Riemann
surfaces of genus $2$}, Trans. Amer. Math. Soc. 144 (1969) 95--113.

\bibitem[BS-a]{s91}Bernhard Schindler, \emph{Jacobische Varietaeten
 hyperelliptischer Kurven und einiger spezieller Kurven vom Geschelcht
 3.}
(Dissertation Erlangen 1991).

\bibitem[BS-b]{s93}Bernhard Schindler, \emph{
Period matrices of hyperelliptic curves}, Manuscripta Math. 78
(1993), no. 4, 369--380.

\bibitem[St]{s01}M. Streit, \emph{Period matrices and representation theory}, Abh.
Math. Sem. Univ. Hamburg 71 (2001), 279--290.

\bibitem[Ta]{t08}Yuuki Tadokoro, \emph{A nontrivial algebraic cycle in
    the Jacobian variety of the Klein quartic}, Math. Z. 260 (2008),
  no. 2, 265--275.

\bibitem[TT]{tretkoff}C. L. Tretkoff and M. D. Tretkoff,  \emph{Combinatorial group
theory, Riemann surfaces and differential equations},
Contributions to group theory,  467--519, Contemp. Math., 33,
Amer. Math. Soc., Providence, RI, 1984.

\bibitem[Yo]{y99}Katsuaki Yoshida, \emph{Klein's surface of genus
    three and associated theta constants},Tsukuba J. Math. 23 (1999),
  no. 2, 383--416.

\end{thebibliography}

\end{document}